\documentclass[12pt]{amsart}
\usepackage{amsmath,amssymb,amsfonts, amscd}
\newtheorem{thm}{Theorem}[section]

\newtheorem{prop}[thm]{Proposition}

\newtheorem{rem}[thm]{Remark}

\newcommand{\iso}{\cong}

\newcommand{\hra}{\hookrightarrow}

\newcommand{\C}{\mathbb C}
\newcommand{\Z}{\mathbb Z}

\newcommand{\pp}{\mathbb P}
\newcommand{\Aut}{\operatorname{Aut}}

\newcommand{\Pic}{\operatorname{Pic}}
\newcommand{\Id}{\operatorname{Id}}
\newcommand{\OO}{\mathcal O}
\newcommand{\OS}{\OO_S}
\newcommand{\fie}{\varphi}

\newcommand{\ga}{\gamma}
\newcommand{\Ga}{\Gamma}
\newcommand{\si}{\sigma}
\newcommand{\Si}{\Sigma}

\newcommand{\bY}{\bar{Y}}
\newcommand{\bT}{\bar{T}}

\newcommand{\tY}{\tilde{Y}}

\newcommand{\tX}{\tilde{X}}
\newcommand{\tC}{\tilde{C}}

\newcommand{\btimes}{\,{\scriptstyle \boxtimes}\,}

\newcommand{\inv}{^{-1}}

\numberwithin{equation}{section}
\title{Regular canonical covers}
\author{Ciro Ciliberto, Rita Pardini, Francesca Tovena}
\date{}
\begin{document}

\begin{abstract}
We construct three sequences of regular surfaces of general type with unbounded numerical
invariants whose canonical map is 2-to-1 onto a canonically embedded surface. Only sporadic
examples of surfaces with these properties were previously known.

\end{abstract}
\maketitle
\section{Introduction}
Let $X$ be a smooth surface of general type and let
$\fie\colon X\to
\Sigma\subseteq {\bf P}^{p_g(X)-1}$ be the canonical map of $X$,  where
$\Sigma$ is the image of $\fie$. Suppose that
$\Sigma$ is a surface and that $\fie$ has degree $d\ge 2$. Let
$\epsilon\colon S\to \Sigma$ be a  desingularization  of
$\Sigma$. A classical result (cf. 
\cite{beauville}, and also \cite{bab}, 
\cite{babbage})
says that either $p_g(S)=0$ or $S$ is of general type and $\epsilon\colon
S\to
\Sigma$ is the canonical map of $S$. In the latter case we have a
dominant rational map $X\to S$ of degree $d$, a so-called 
{\it good canonical cover} of degree $d$ (see \cite{cirifra}).\par

A basic construction of sequences of good canonical covers with unbounded invariants is
due to Beauville (see \cite{Ca}, and
\cite{MP}). In \cite{cirifra} (see also \S 2) we have presented a more general construction
which produces new infinite series of good canonical covers $X\to S$ of degree 2.  All
the examples obtained by this method  have $q(X)\ge 2$ and,  although some  examples of
canonical covers $X\to S$ with $X$ regular can be found in the literature (see \cite{zagier},
\cite{beauville} Proposition 3.6, \cite{babbage} Theorem 3.5, \cite{enriques},
\cite{supcan}), no
infinite series of such examples has been presented so far.

In this note we fill up this gap. Our idea, roughly speaking, is to take a good
canonical cover $X\to S$ and find 
a group $G$ acting on both $X$ and $S$ in such a way that $X/G\to
S/G$ turns out to be again a good canonical cover and moreover $q(X/G)=0$. This idea 
works nicely: for each infinite series of good canonical covers produced in \cite{cirifra} we
are able to find a $\Z_3-$action which does the job (see \S 3).
For one series
of the above examples, we prove in addition that $X$ is simply connected (see \S 4).
\bigskip

\paragraph{\bf Acknowledgments:} The second author wishes to thank Daniel Naie for
interesting conversations on the subject of this paper.

This research was carried out in the framework of the EU Research Training Network
EAGER. The authors are
members of G.N.S.A.G.A. of C.N.R.
\bigskip

\paragraph{\bf Notation:} All varieties are  projective varieties over the field of ccomplex
numbers.  We write
$\omega:=e^{2\pi i/3}$.

A {\em map} is a rational map and a {\em morphism} is a regular map.  We do not
distinguish between Cartier divisors and line bundles and use the additive
and
multiplicative notation interchangeably.  If $X:=X_1\times X_2$
is a product of varieties with projections  $p_i\colon X\to X_i$ and
$L_i$ is a line bundle on $X_i$, $i=1,2$, then $L_1\btimes L_2$ denotes the line bundle
$p_1^*L_1\otimes p_2^*L_2$. 

The remaining   notation is standard in  algebraic geometry; we
just recall here the notation for the invariants of a smooth  surface $S$: $K_S$
is the {\em canonical class}, $p_g(S)=h^0(S,K_S)$ the {\em geometric
genus}, $q(S)=h^1(S,\OS)$ the {\em irregularity}, and $\chi(S)=1+p_g(S)-q(S)$ the {\em Euler
characteristic}.

\section{Good canonical covers and generating pairs}
In this section we recall the definitions and the results that we need from \cite{cirifra} in
the form most suitable for our purposes.

A {\em canonical cover} of degree $d\ge 2$ is a generically finite rational map $X\to S$ of
degree $d$  between smooth surfaces of general type such that $p_g(X)=p_g(S)$. If in addition
the canonical map of $S$ is birational onto its image, then we say that the canonical cover
is {\em good}.  By Proposition 4.1 of \cite{beauville}, one has $d\le 3$ for large values
of the numerical invariants of $X$.

In \cite{cirifra}, generalizing a construction of Beauville, it is shown how to  construct an
infinite series of good canonical covers of degree 2 starting with  a so-called  {\em good
generating pair}. A  good generating pair is  a pair 
$(h\colon V\to W, L)$, where $h$ is a finite morphism  of degree 2 between surfaces and $L$ is
a line bundle on
$W$ that satisfy  a certain  set  of conditions (see \cite{cirifra},
Definition 2.4 for the precise definition).  In particular the surface $V$ is smooth, the
involution $\iota$ acts freely in codimension 1 and $W$ is normal, hence $h\colon V\to
 W=V/\iota $ is the quotient map and the singularities of $W$ are $A_1$ points, which are the
images of the fixed points of $\iota$. The general curve $C$ of $|L|$ is a smooth connected
non hyperelliptic  curve, whose genus $g$ is called the genus of the generating pair.

 The main point of \S 2 of \cite{cirifra} is that a  sequence of
good  canonical covers can be obtained from a good generating pair in the following way.
Denote by
${\mathcal L}(n)$ the line bundle $h^*L\btimes\OO_{\pp^1}(n)$ on $V\times \pp^1$. For $n\ge
3$ we take 
$X\subset V\times \pp^1$  to be a smooth element of
$|{\mathcal L}(n)|$. By the properties of good generating pairs,  $X$ is mapped to itself by
the involution
$\iota\times Id$ and the quotient surface $\Sigma\subset W\times \pp^1$ is an element of the
linear system $|L\btimes \OO_{\pp^1}(n)|$ with
$A_1$ singularities. If we denote by
$S$ the minimal resolution of $\Si$, then  the induced rational map $X\to S$ is a
good canonical cover of degree 2. The composite map $S\to W\times \pp^1\to\pp^1$ is a
fibration whose general   fibre $F$ is a non hyperelliptic curve of genus $g$ (the map
$S\to W\times \pp^1\to W$ identifies $F$ with a curve of $|L|$). The surface
$S$ is regular, while the surface $X$ has irregularity $g-1\ge 2$. The 
invariants $\chi(X)$ and $K^2_X$ grow linearly with $n$, hence one  obtains good canonical
covers with arbitrarily large invariants.
\bigskip

We describe now briefly the three   examples of generating pairs we know of. A more detailed
description is contained in \S3 of \cite{cirifra}. 
\smallskip

\noindent{\bf Generating pair I:} Let $C$ be a smooth projective curve of genus $2$. We set
$V:=\Pic^1(C)$, we let $\iota$ be the involution of $V$ defined by $N\mapsto K_C\otimes N\inv$
and we let
$h\colon V\to W:=V/\iota$ be the quotient map. The involution
$\iota$ has 16 fixed points.  The map
$C\hra V$ defined by $Q\mapsto \OO_C(Q)$ identifies $C$ with the set $\{N\in
\Pic^1(C)|h^0(N)>0\}$, and the  line bundle
$L$ on $W$ is determined by the condition that
$h^*L=\OO_V(2C)$ (up to non canonical isomorphism, $V$ is an
 abelian surface with an  irreducible principal polarization and $W$ is its Kummer surface).

 This  pair has genus $g=3$, and 
 the numerical invariants of the corresponding good canonical covers are the following:
$$ q(X)=2;\quad  p_g(X)=4n-3;\quad K^2_X=24n-32$$
$$ q(S)=0;\quad p_g(S)=4n-3;\quad K^2_S=12n-16.$$
\smallskip

\noindent{\bf Generating pair II:} Let $C$ be a smooth curve of genus 2 and, as in the
description of the generating pair I, consider the natural embedding $C\hra \Pic^1(C)$. We set
$M:=\OO_{\Pic^1(C)}(C)$, we let
$D$ be a smooth divisor of $|2M|$ and take $\pi\colon V\to\Pic^1(C)$ to be the double cover of
$\Pic^1(C)$ branched on $D$ and  such that $\pi_*\OO_V=\OO_{\Pic^1(C)}\oplus M\inv$. The
surface
$V$ is smooth and the involution of
$\Pic^1(C)$ defined by $N\mapsto K_C\otimes N\inv$ can be lifted to an involution $\iota$ of
$V$ with 20 fixed points. We take $h\colon V\to W:=V/\iota$ to be the quotient map and we set
$L:=K_W$.  This  pair has genus $g=3$, and 
 the numerical invariants of the corresponding good canonical covers are the following:
$$ q(X)=2;\quad p_g(X)=5n-3;\quad K^2_X=32n-32$$
$$ q(S)=0;\quad p_g(S)=5n-3;\quad K^2_S=16n-16.$$
\smallskip

\noindent{\bf Generating pair III:} Let $C$ be   a smooth non hyperelliptic projective  
curve of genus 3. We set $V:=S^2C$,  we let $\iota$ be the involution of $V$
defined by  $Q+R\mapsto K_C(-Q-R)$,  we let $h\colon V\to W:=V/\iota$ be the quotient map.
 The involution $\iota$ has 28 fixed points. We take $L:=K_W$.  This  pair has genus $g=4$,
and 
 the numerical invariants of the corresponding good canonical covers are the following:
$$ q(X)=3;\quad p_g(X)=7n-4;\quad K^2_X=48n-48$$
$$ q(S)=0;\quad p_g(S)=7n-4;\quad K^2_S=24n-24.$$
\smallskip

All the examples constructed in the next section involve  taking quotients of  surfaces
by a $\Z_3-$action with isolated fixed points. Next we give  the formulas that we  need to
compute the invariants of the minimal resolution of such a quotient. We recall that by
Cartan's lemma  the representation of $\Z_3$ on the tangent space at an isolated fixed point
$Q$ is the sum of two nontrivial characters. The  image
of $Q$ in the quotient surface is a canonical  singularity of type $A_2$ if the two 
characters are different, and it is a $\frac{1}{3}(1,1)$ singularity, i.e. a rational singularity   solved by   a smooth
rational curve of self intersection
$-3$, if the two
characters are equal.
\begin{prop}\label{invariants}
Let  $X$ be a projective surface with canonical singularities. Assume that
$\Z_3$ acts on
$X$ with finitely many fixed points, which are all smooth for $X$,  and denote by
$Y$  the minimal resolution of
$X/\Z_3$.  If $\alpha$ (resp. $\beta$)  is the number of fixed points
$Q$ of
$\Z_3$ on
$X$ such that the characters of the representation of $\Z_3$ on $T_QX$ are equal (resp.
distinct), then:
$$K^2_X=3K^2_Y+\alpha;\quad \chi(X)=3\chi(Y)-\alpha/3-2\beta/3.$$
\end{prop}
\begin{proof}
First of all it is easy to show that  $\Z_3$ acts on the minimal desingularization of $X$ and
that the quotient  of the minimal desingularization is a partial resolution of $X/\Z_3$. In
addition the invariants of $X$ are not affected by the resolution, since $X$ has canonical
singularities. Hence we may assume that $X$ is smooth.

 Let
$Q\in X$ be a fixed point of
$\Z_3$.  Assume that $\Z_3$ acts on $T_QX$ with equal characters. If we blow up $Q$, then
$\Z_3$ acts on the blown up surface fixing the exceptional curve pointwise.
If instead $\Z_3$ acts on $T_QX$ with distinct characters, then the action on the blown up
surface has two isolated fixed points $Q_1$, $Q_2$  on the exceptional curve, corresponding
to the $\Z_3-$eigenspaces of $T_QX$. A local computation shows that $\Z_3$ acts with
equal characters on the tangent space at $Q_1$ and $Q_2$. Summing up, it is possible to blow
up $X$ in such a way that $\Z_3$ acts on the blow up $\tX$ and the fixed locus of $\Z_3$ on
$\tX$  is a divisor. It follows that  $\tY:=\tX/\Z_3$ is smooth and the quotient map
$\tX\to
\tY$ is a flat morphism. There is a commutative diagram:
\[
\begin{CD}
\tX @>>> X\\
@VVV    @VVV\\
\tY @>>> X/\Z_3
\end{CD}
\]
where the map $\tX\to X$ is the blow up at $\alpha+3\beta$ points and $\tY\to X/ \Z_3$
is the composition of the resolution $Y\to X/\Z_3$ with the map $\tY\to Y$ obtained by
blowing  up the $\beta$  intersection points of  the  pairs of $-2$ curves in the resolution
of the
$A_2-$singularities of $X/\Z_3$.
 The statement now follows by using  the
standard formulas for cyclic   covers of smooth surfaces to relate  the invariants of $\tX$
and
$\tY$.
\end{proof}
\section{  Construction of the regular good canonical covers}
In this section we construct infinite series of canonical covers $Y\to
T$ of degree 2 such that $q(Y)=0$. We use the notation and the terminology introduced in \S 2.

Our strategy is the following. Let $X\to \Si$ be a finite degree 2 map of surfaces, with $X$
smooth and $\Si$ normal, denote by $S$ the minimal desingularization of $\Si$ and assume that
the induced map $X\to S$ is a good canonical cover. Let $G$ be a subgroup of $\Aut(X)$ that
does not contain the canonical involution $\si$ of $X$. Since $\si$ is
in the center of $\Aut(X)$, it follows that $G$ acts on $\Si=X/\si$.
Set $\bY:=X/G$ and $\bT=\Si/G$ and denote by $\bar{\pi}\colon \bY\to \bT$ the induced map,
which is a finite morphism of degree 2. The surfaces
$\bY$ and $\bT$ have rational singularities.  If $Y$ is the minimal desingularization of
$\bY$ and  $T$ is the minimal desingularization of $\bT$, then there  are   natural
identifications (cf. \cite{steenbrink}, Lemma 1.8 and  1.11):  
$H^i(Y, \Omega_Y^j)\iso H^i(X,\Omega_X^j)^G$ and
$H^i(T, \Omega^j_T)\iso H^i(\Si,\Omega^j_{\Si})^G$, for $0\le i,j\le 2$.
It
follows  that $p_g(Y)=p_g(T)$, namely the canonical map of $Y$ factors through the induced
map 
$Y\to T$.  If in addition  the canonical map of $T$ is
birational and
$H^0(X,\Omega^1_X)^G=\{0\}$, then we have the required example.

We apply this idea to  canonical covers obtained from the three generating pairs of
 \S 2, using $\Z_3-$actions  and giving an example for each one.  
One  can obtain similar examples by using
 other groups.  It is also possible to construct by the same method infinite series
of good canonical covers
$Y\to T$ such that
$q(Y)=1$.
\smallskip

\noindent{\bf Example 1:} Here we take $\Z_3-$quotients of  good canonical covers obtained
from the generating pair I (cf. \S 2).

Assume that the curve $C$ of genus 2 has  a 
$\Z_3-$action such that $C/\Z_3$ is rational (such a curve can be contructed using for
instance  Proposition 2.1 of
\cite{ritaabel}).  By the
Hurwitz formula,  $\Z_3$ has  $4$ fixed
points on $C$. Denote by
$\ga$ the hyperelliptic involution of
$C$ and by $\psi\colon C\to C/\ga= \pp^1$ the canonical map. Since
$\ga$ commutes with all the automorphisms of $C$, there is an induced action
of
$\Z_3$ on $\pp^1$ that preserves the branch locus $B$ of $\psi$. 
Since $B$ consists of $6$ points, it follows that it is the union of two orbits of the
$\Z_3$-action on $\pp^1$ and therefore the fixed points $P_0$, $P_1$ of 
$\Z_3$ on $\pp^1$ do not belong to  $B$.
 Let $x_0, x_1$ be homogeneous
coordinates on $\pp^1$ such that $P_0=(1:0)$ and $P_1=(0:1)$ and let
$\tau_i=\psi^*x_i$, $i=0,1$.  Then
$\tau_0$,
$\tau_1$ is a basis of
$H^0(C,\omega_C)$ and we denote by $\xi$ the generator of $\Z_3$ such that
$\xi^*\tau_0=\omega \tau_0$ and $\xi^*\tau_1=\omega^2\tau_1$, where $\omega:= e^{2\pi
i/3}$. We
write $Q_1+Q_2$ for the divisor of zeros of $\tau_0$ on $C$ and $Q_3+Q_4$ for the
divisor of zeros of $\tau_1$. The points $Q_1\ldots Q_4$ are distinct and they
are the fixed points of  $\Z_3$. Clearly, $\ga$ maps $Q_1$ to
$Q_2$ and $Q_3$ to $Q_4$.

The $\Z_3$-action on $C$ extends to  a
$\Z_3$-action on $\Pic^1(C)$ compatible with the natural embedding $C\hra \Pic^1(C)$. Set
$M:=\OO_{\Pic^1(C)}(C)$. In the notation of \S 2, we have $M^{\otimes 2}=h^*L$.  The line
bundle
$M$  can be identified with a subsheaf of the constant sheaf
$\C(\Pic^1(C))$,  hence the    action of
$\Z_3$ on  $\C(\Pic^1(C))$  restricts to a linearization of $M$ such that  
$\Z_3$ acts trivially on the $1-$dimensional vector space $H^0(\Pic^1(C), M)$. We consider on
$M^{\otimes 2}$ the $\Z_3-$linearization induced by that of $M$.

The above description of the $\Z_3-$action on $H^0(C, \omega_C)$ shows
that 
$\Z_3$ acts trivially on
$\omega_{\Pic^1(C)}=\OO_{\Pic^1(C)}$. Consider the residue sequence:
$$0\to M\to M^{\otimes 2} \to \omega_C^{\otimes 2}\to 0.$$ 
 It is easy to check that
  the maps in the sequence  are equivariant with respect to the chosen
actions on
$M$ and
$M^{\otimes 2}$ and with respect to the natural action on
$\omega_C^{\otimes 2}$. The space $H^0(C,\omega_C^{\otimes 2})$ is spanned by
$\tau_0^2,\tau_1^2,
\tau_0\tau_1$ and one has: $\xi^*\tau_0\tau_1=\tau_0\tau_1$, 
$\xi^*\tau_0^2=\omega^2\tau_0^2$ and $\xi^*\tau_1^2=\omega \tau_1^2$.  Since the restriction  map
on global sections
$H^0(\Pic^1(C),M^{\otimes 2})\to H^0(C,\omega_C^{\otimes 2})$ is surjective  by Kodaira
vanishing, one  can  find eigenvectors 
$f_1,f_2, f_3\in H^0(\Pic^1(C), M^{\otimes 2})$ that map to $\tau_0\tau_1$, $\tau_0^2$, $\tau_1^2$,
respectively. It follows that a basis of eigenvectors for 
$H^0(\Pic^1(C), M^{\otimes 2})$ is given  by
$f_0, f_1, f_2, f_3$, where $f_0$ is the square of a nonzero section of
$M$, and thus it   vanishes on
$C$ of order $2$. 
The base scheme of the pencil of $\Pic^1(C)$ spanned by $f_0$ and $f_1$ is supported at the
points
$Q_1, Q_2, Q_3, Q_4$ of $C$ and it has length $2$ at each of these points. It
follows that the general curve of this pencil  is smooth  by Bertini's theorem.

Consider the $\Z_3-$action on $\pp^1$ defined by
$\xi(x_0:x_1)=(\omega x_0: \omega^2 x_1)$ and let $\Z_3$ act diagonally on
$\Pic^1(C)\times
\pp^1$ by
$\xi(a,b)=(\xi a, \xi b)$. We remark that  this action has isolated fixed points.  If we
extend the
$\Z_3$-action on $\Pic^1(C)\times \pp^1$ to ${\mathcal L}(3k):=M^{\otimes 2}\btimes \OO_{\pp^1}(3k)$ by
using the
$\Z_3$-action that we have already defined on $M$ and the natural pull-back structure on
$\OO_{\pp^1}(3k)$, then  the following is a basis of the $\Z_3$-invariant
subspace of $H^0(\Pic^1(C)\times \pp^1,{\mathcal L}(3k))$:

$$x_0^{3i}x_1^{3k-3i}f_0,\ \ x_0^{3i}x_1^{3k-3i}f_1,\ \ 
x_0^{2+3j}x_1^{3k-3j-2}f_2,\ \  x_0^{1+3j}x_1^{3k-3j-1}f_3,$$
$$0\le i\le k,\ \  0\le j\le k-1$$
Denote by $|X|$ the corresponding linear system on $\Pic^1(C)\times \pp^1$. The base
set of  $|X|$ consists of the $8$ points  
$\{(Q_i, (1: 0)),(Q_i,(0:1))|\  i=1\ldots 4\}$,  which are
all fixed points of $\Z_3$. The general surface of $|X|$ contains no other
fixed point and one can check that it is smooth using   Bertini's
theorem. We wish to determine the characters of the representation of
$\Z_3$ on the tangent space  to a general $X$ at the fixed points of
$|X|$.  We remark 
 first of all that, since  the tangent space  to $\Pic^1(C)$ is dual to
 $H^0(\omega_C)$, the $\Z_3-$action on it is the sum of the two nontrivial
characters.
Let $T'$
be the tangent space to $X$ at $(Q_1,(0:1))$. Since
$f_3$ does not vanish at $Q_1$, for a general choice of $X\in |X|$ the
projection $X\to \Pic^1(C)$ induces an isomorphism of $T'$ with the tangent space to
$\Pic^1(C)$. Since $X\to \Pic^1(C)$ is compatible with the $\Z_3-$actions on $X$ and
$\Pic^1(C)$, it follows that $\xi$ acts on  $T'$ with eigenvalues $\omega$ and $\omega^2$.
 Now let $T''$ be the tangent space to $X$ at the point $(Q_1, (1:0))$. 
We have  seen that  the tangent space at $Q_1$  to  the curves of the pencil of
$\Pic^1(C)$ spanned by
$f_0$,
$f_1$ is fixed and transversal to $C$.  We denote  by $T_1$ this space, which is clearly 
$\Z_3-$invariant.   Since
$f_2$ vanishes at
$Q_1$, it follows that $T''$ is the direct sum of $T_1\times\{0\}$  and of the tangent
space
 to $Q_1\times \pp^1$ at $(Q_1, (1:0))$, on which $\xi$ acts  with eigenvalue  $\omega$.
We observe that $\tau_1$ is a generator of $\omega_C$ in a neighbourhood of
$Q_1$, hence     $\xi$ acts by pull-back on  the cotangent space to $C$
at
$Q_1$ as the multiplication by  $\omega^2$. Thus  $\xi$ acts by
differentiation on the tangent space $T_2$ to $C$ at $Q_1$    as the multiplication
by $\omega^2$. Since  the tangent space to $\Pic^1(C)$ at $Q_1$ is the direct sum of $T_1$ and
$T_2$, it follows that
$\xi$ acts on $T_1$ as the multiplication by $\omega$. Summing up, we have shown
that the action of $\xi$ on $T''$ is the multiplication by $\omega$.    The
representations of
$\Z_3$ at the remaining six points can be studied in the same way, and therefore the
characters of the representation of $\Z_3$ 
on the tangent space are distinct at the points  $(Q_1, (0:1)),
(Q_2,(0:1))$, $(Q_3, (1:0)), (Q_4,(1:0))$ and they are equal at 
the points $(Q_1, (1:0)), (Q_2, (1:0)), (Q_3,
(0:1)),(Q_4,(0:1))$.  Let $Y$ be the minimal desingularization of $X/\Z_3$. Then by
Proposition \ref{invariants} we have:
$$K^2_X-4=3K^2_{Y};\quad \chi(X)=3\chi(Y)-\frac{4}{3}-\frac{8}{3}.$$
In addition, by  the Lefschetz Theorem on hyperplane sections, the
Albanese variety of $X$ is $\Pic^0(C)$ and thus $q(Y)=0$, since $\Z_3$ acts on $\Pic^0(C)$
with finitely many fixed points.
Using the formulas for the  invariants of $X$ (cf. \S 2), we get: 
$$K^2_Y=24k-12;\quad q(Y)=0;\quad p_g(Y)= 4k-1.$$
As explained in \S 2,  the canonical map of $X$ has degree 2 and the corresponding 
involution $\si$  of $X$ is induced by the involution 
$\iota\times Id$ of
$\Pic^1(C)\times \pp^1$. The quotient surface  
$\Si:=X/\si$ is a normal regular surface with  $48k$ singular points of type $A_1$ that are
the images of the fixed points of $\iota\times Id$ on $X$. The
$\Z_3$-action on
$X$ commutes with
$\si$ and therefore there is an induced $\Z_3$-action on $\Si$. 
The involution
$\iota$ exchanges
$Q_1$ with $Q_2$ and $Q_3$ with $Q_4$, and thus there are $4$ fixed points of
$\Z_3$ on $\Si$, that are smooth points of $\Si$. The eigenvalues of the
action of $\xi$ on the tangent space to $\Si$ are equal at two of these points
and they are different at the remaining two. 
If $T$ is the minimal desingularization of $\Si/\Z_3$, then by Proposition \ref{invariants} we
have:
$$K^2_{\Si}-2=3K^2_T; \quad \chi(\Si)=3\chi(T)-2/3-4/3.$$
The invariants of $\Si$ are computed in \S 2, and we get:
$$K^2_T=12k-6;\quad q(T)=0;\quad p_g(T)=4k-1 .$$
Hence $p_g(Y)=p_g(T)$.
Finally, we are going to show that the canonical map of $T$ is birational
for $k\ge 2$. It is easy to check that the non hyperellyptic  genus 3 fibration $\Si\to \pp^1$
induced by the second projection $W\times \pp^1\to \pp^1$ descends to a fibration $\Si/\Z_3\to
\pp^1/\Z_3$ with the same general fibre $F$, which in turn pulls back to a fibration $p\colon
T\to\pp^1$. For
$k\ge 2$, we have $p_g(T)\ge 7$, hence the canonical map of $T$ separates the fibres of $p$.
Thus, in order to show that the canonical map of $T$ has degree 1 it is enough to show that
the restriction of $|K_T|$
 to the general fibre of $p$ is the complete canonical system of the fibre.

The surface $\Si\subset W\times \pp^1$ is an element of the linear system
$|L\btimes\OO_{\pp^1}(3k)|$. The sheaf $L\btimes\OO_{\pp^1}(3k)$ inherits from $M^{\otimes
2}\btimes
\OO_{\pp^1}(3k)$ a linearization such that $\Si$ is defined by a $\Z_3-$invariant section.
The residue sequence:
$$0\to K_{W\times\pp^1}\to K_{W\times\pp^1}(\Si)=L\btimes\OO_{\pp^1}(3k-2)\to K_{\Si}\to 0$$
is $\Z_3-$equivariant  with respect to the natural linearization of $K_{W\times\pp^1}$ and
to  the chosen linearization of
$L\btimes\OO_{\pp^1}(3k-2)$.

For $i=0\dots 3$, let $g_i\in H^0(W, L)$ be
such that $h^*g_i=f_i$ (recall that we have $H^0(V, M^{\otimes 2})=h^*H^0(W,L)$). 
  A
 basis of the invariant subspace $H$ of $H^0(W\times \pp^1,L\btimes 
\OO_{\pp^1}(3k-2))$ is the following:
$$x_0^{2+3i}x_1^{3k-4-3i}g_0,
x_0^{2+3i}x_1^{3k-4-3i}g_1, x_0^{1+3m}x_1^{3k-3-3m}g_2,
 x_0^{3m}x_1^{3k-2-3m}g_3$$
$$0\le i\le k-2,\,\, 0\le  m \le k-1.$$ 
Consider a general fibre $F$ with $(x_0:x_1)=(\lambda_0: \lambda_1)$. Identify $F$  with
the corresponding curve $C$ of $|L|$ on $W$. Then the above sections, restricted to $F$,
correspond to the restrictions of $g_0\dots g_3$ to $C$. Hence,  the restriction map $H\to
H^0(F,K_F)$ is surjective.  Since the  elements of
$H$ give  sections of the canonical bundle  of $T$ and the fibrations
$\Si\to\pp^1$ and $p\colon T\to\pp^1$ have the same general fibre, it follows  that the
restriction of $|K_T|$ to the general fibre of $p$ is also complete.
\bigskip

\noindent {\bf Example 2:}  Here we take $\Z_3-$quotients of  good canonical covers obtained
from the generating pair II (cf. \S 2).

We take   a curve $C$ of genus 2 with a 
 $\Z_3-$action with rational quotient as in Example 1  and we use the notation
introduced there. 

We take $\pi\colon V\to \Pic^1(C)$ to  be a double cover  such that
$\pi_*\OO_V=\OO_{\Pic^1(C)}\oplus M\inv$,  branched on a general 
 curve $D$ of $|M^{\otimes 2}|$ invariant under the $\Z_3-$action.  We have seen in Example 1
that such curves are a
$1-$dimensional linear system with simple  base points
$Q_i$, $i=1\dots 4$. Using the chosen $\Z_3-$action on $M$ it is possible to lift the
automorphism
$\xi$ of
$\Pic^1(C)$ to
$V$, hence there is an automorphism $\zeta$ of $V$ of order $6$ such that $\zeta^2$
lifts $\xi$ and $\zeta^3$ is the involution associated to the double cover $\pi\colon
V\to \Pic^1(C)$.  The fixed points of $\zeta^2$ are isolated, since  the fixed
points of
$\xi$ on $\Pic^1(C)$ are isolated, and they are of two types: the 10 inverse images
$P_1\dots P_{10}$  of the fixed points of $\xi$ not lying on $B$ and the points
$R_i:=\pi\inv Q_i$,
$i=1\dots 4$. Since the differential of $\pi$ is nonsingular at $P_1\dots
P_{10}$, $\xi=\zeta^2$ acts on the tangent space to $V$ at $P_1\dots P_{10}$ with
eigenvalues $\omega$, $\omega^2$. The points $R_1\dots R_4$ are fixed points of
$h$.
 Now consider the minimal resolution $Z$ of the quotient surface $V/\xi$.  Recall
that
$V$ has invariants $q(V)=p_g(V)=2$ and $K^2_V=4$. By Proposition
\ref{invariants}, one has $4=K^2_V=\alpha+3K^2_Z$ and 
$1=\chi(V)=3\chi(Z)-(\alpha+2\beta)/3$, where $\alpha$ (resp. $\beta$) is the number of
isolated fixed points of $\xi$ on $V$ such that the eigenvalues of the action of $\xi$ on the
corresponding tangent space are equal (resp. distinct). By the above discussion, we have
$\beta\ge 10$ and $\alpha+\beta=14$, hence  the only possibility is
$\alpha=4$,
$\beta=10$. Thus we have
$\chi(Z)=3$,
$p_g(Z)=2$, $q(Z)=0$, $K^2_Z=0$. 
In particular, $p_g(Z)=p_g(V)$, namely the action of $\Z_3$ on $H^0(V, K_V)$ is trivial. 
We let
$\xi$ act on
$\pp^1$ by $\xi(x_0,x_1)=(\omega x_0, \omega^2 x_1)$ and we consider the diagonal action on 
$V\times \pp^1$. If $f_0$, $f_1$ is
a basis of $H^0(V, K_V)$ and we consider the natural action of $\xi$ on ${\mathcal
L}(3k)=K_V\btimes\OO_{\pp^1}(3k)$, then the following is a basis of the subspace  of
invariant sections of
${\mathcal L}(3k)$ with respect to the induced linearization:
$$x_0^{3i}x_1^{3k-3i}f_0,\quad x_0^{3i}x_1^{3k-3i}f_1;\quad i=0\dots k.$$
 It is easy to check that the base locus of the corresponding linear system
$|X|\subset |{\mathcal L}(3k)|$ is the union of the curves $\{R_i\}\times \pp^1$,
$i=1\dots 4$, and that the general $X$ is smooth along the base locus. Hence the
general $X$ is smooth by Bertini's theorem.  The fixed
points of the action of
$\xi$ on the general $X$ are the points $(R_i,(0:1))$ and $(R_i,(1:0))$,
$i=1\dots 4$. Arguing as in the previous case, it is not difficult to check that the
irreducible characters of the representation of $\Z_3$ on the tangent space to $X$ are equal
at four of these points
 and they are different  at the remaining four. 
Denote by $Y$ the minimal desingularization of $X/\Z_3$.  By the formulas for the invariants
of $X$ (cf. \S 2) and by Proposition \ref{invariants}, we have: 
$$\chi(Y)=5k;\quad K^2_Y=32k-12.$$
In addition, we have $q(Y)=0$, since the Albanese variety of $Y$ is $\Pic^0(C)$ and
$\Z_3$ acts on $\Pic^0(C)$ with isolated fixed points. Thus  we have $p_g(Y)=5k-1$.
As explained in \S 2, the canonical map of $X$ has degree 2 and the canonical
involution $\si$ of $X$ is the restriction of $\iota\times Id$. 
The action of $\xi$ on $X$ descends to $\Si=X/\si$ and  $\xi$ has 4 fixed points occurring
at smooth points $\Si$.  The action of $\Z_3$ on the tangent space to $\Si$ has distinct
 characters at two of these points  and equal characters at the remaining two. The
invariants of the minimal resolution $T$ of $\Si/\Z_3$, which  can be computed as before, 
are:
$$q(T)=0;\quad p_g(T)=5k-1;\quad K^2_T=16k-6.$$
Hence we have $p_g(Y)=p_g(T)$.

One can
show that for $k\ge 1$ the canonical map of $T$ is birational, and thus $Y\to T$ is a
canonical cover, by using an  argument similar to that of Example 1. In order to do so, take a
basis $g_0\dots g_4$ of $H^0(W,2K_W)$ consisting of $\Z_3-$eigenvectors. The $\Z_3-$invariant
subspace  $H$ of the space  $H^0(W\times\pp^1,2K_W\btimes \OO_{\pp^1}(3k)$ is spanned by
products of the
$g_i$ with suitable monomials in $x_0,x_1$ and each of the $g_i$ occurs in at least one such
product. Consider a general fibre $F$ of the map $\Si\subset  W\times \pp^1\to \pp^1$, with
$(x_0:x_1)=(\lambda_0:
\lambda_1)$. Identify
$F$  with the corresponding curve $C$ of $|K_W|$ on $W$. Then  the elements of the  above
mentioned  basis of
$H$, restricted to
$F$, correspond to the restrictions of $g_0\dots g_4$ to $C$. Hence,  the restriction map
$H\to H^0(F,K_F)$ is surjective.  Since the  elements of
$H$ give  sections of the canonical bundle  of $T$ and the fibrations
$\Si\to\pp^1$ and $p\colon T\to\pp^1/\Z_3=\pp^1$ have the same general fibre, it follows  that
the restriction of $|K_T|$ to the general fibre of $p$ is also complete.
\bigskip

\noindent{\bf Example 3:}  Here we take $\Z_3-$quotients of  good canonical covers obtained
from the generating pair III (cf. \S 2).

 Consider $\pp^2$ with homogeneous coordinates
$z_0,z_1,z_2$ and  let $C\subset \pp^2$ be the curve defined by $z_1z_2^3=f(z_0,
z_1)$ where
$f$ is a homogeneous polynomial of degree $4$ with distinct roots and such that
$f(1,0)\ne 0$. $C$ is a smooth non hyperelliptic
curve of genus 3 and $(z_0:z_1:z_2)\mapsto (z_0:z_1:\omega z_2)$ defines an
automorphism $\xi$ of $C$  of order 3. The fixed points of $\xi$ are $Q_0:=
(0:0:1)$ and the points $Q_1,\dots Q_4$ defined by $z_2=f(z_0,z_1)=0$. 
We have $H^0(C,\omega_C)=H^0(C,\OO_C(1))$ and $\xi$ acts on $1-$forms as
follows: 
$$\xi^*z_0=\omega z_0,\quad
\xi^*z_1=\omega z_1,\quad \xi^*z_2=\omega^2 z_2.$$
The induced $\Z_3-$action on $V=S^2C$ clearly commutes with the involution $\iota$ of $V$.
In order to write down the action of $\Z_3$ on  $H^0(V, K_V)$, we observe that there
is a natural isomorphism   $\bigwedge^2
H^0(C,\omega_C)\cong H^0(V,K_V)$. If we denote by $f_i$ the element of $H^0(V,K_V)$
corresponding to $z_j\wedge z_k$, with $(ijk)$ a permutation of $(012)$, then we have:
$$\xi^*f_0=f_0;\quad \xi^*f_1=f_1;\quad \xi^*f_2=\omega^2 f_2.$$ 
In addition, 
  given a vector space $U$ of
dimension $3$, the alternation  map $U\otimes\bigwedge^2U\to \bigwedge^3U$
gives a natural isomorphism $\bigwedge ^2 U\cong U^*\otimes (\bigwedge^3
U)^*$.  For $U=H^0(C,\omega_C)$ this gives a projective duality between 
$\pp^2=\pp(H^0(C,\omega_C)^*)$ and
${\pp^2}^*=\pp(H^0(V,K_V)^*)$ such that $z_0, z_1, z_2$ and $f_0, f_1, f_2$ are dual
homogeneous coordinates. Hence  the canonical map
$\varphi\colon V\to
{\pp^2}^*$ maps
$P+Q\in V$ to the line $<P,Q>$.  
In particular, the zero locus of $f_0$ on $V$ is the set of the $P+Q$
such that the line determined by $P$ and $Q$ goes through $(1\!:\!0\!:\!0)$, the zero locus of
$f_1$ is the set of
$P+Q$  such that the line determined by $P$ and $Q$ goes through
$(0\!:\!1\!:\!0)$ and the zero locus of $f_2$ is the set of
$P+Q$  such that the line determined by $P$ and $Q$ goes through
$(0\!:\!0\!:\!1)$.
As in the previous examples, we
extend the $\Z_3-$action to $V\times \pp^1$ by setting $\xi(P+Q,
(x_0:x_1))=(\xi (P+Q), (\omega x_0: \omega^2 x_1))$. The fixed points of this
action are the points  $(Q_i+Q_j, (1:0))$,
$(Q_i+Q_j, (0:1))$, $0\le i,j\le 4$. For $k\ge 1$,  we consider  on ${\mathcal
L}(3k)=K_V\btimes
\OO_{\pp^1}(3k)$ the $\Z_3-$linearization induced by the natural one on $K_V$ and by the
linearization of $\OO_{\pp^1}(3k)$ as a pull-back.  The invariant subsystem $|X|$ of 
$|{\mathcal L}(3k)|$ is spanned by:
$$x_0^{3i}x_1^{3k-3i}f_0,\ \ x_0^{3i}x_1^{3k-3i}f_1,\ \ 
 x_0^{1+3m}x_1^{3k-3m-1}f_2,$$
$$0\le i\le k,\ \  0\le m\le k-1$$
The base scheme of $|X|$ consists of the $12$ simple points $(Q_i+Q_j,(1:0))$,
$(Q_i+Q_j,(0:1))$, $1\le i, j\le 4$, $i\ne j$, hence the general surface of $|X|$ is smooth
by Bertini's theorem. The differential of $\xi$ on the tangent space to $V$ at $Q_i+Q_j$ $i\ne
j, i,j\ne 0$,  is the multiplication by $\omega$. For a general $X$ the projection $X\to V$
has non singular differential at the points $(Q_i+Q_j,(0:1))$, $i\ne j$, $i,j\ne 0$, hence the
differential of $\xi$ at these points is also the multiplication  by $\omega$.
Arguing as in Example 1, one can also show that the differential of $\xi$ on the 
tangent space to
$X$ at the points $(Q_i+Q_j,(1:0))$ has eigenvalues $\omega$ and $\omega^2$. 

The invariants of the minimal desingularization $Y$ of $X/\Z_3$ can be computed as in the
previous examples and they are the following:
$$K^2_Y=48k-18;\quad \chi(Y)=7k.$$
In addition, there are natural identifications  $H^0(X,\Omega^1_X)\cong H^0(V,
\Omega^1_V)\cong H^0(C,\omega_C)$, compatible with the various $\Z_3-$actions. It follows that
$H^0(V,\Omega^1_V)^{\Z_3}=\{0\}$, and:
$$q(Y)=0;\quad p_g(Y)=7k-1.$$
As explained in \S 2,  the canonical map of $X$ has degree 2 and the corresponding 
involution $\si$  of $X$ is induced by the involution 
$\iota\times Id$ of
$V\times \pp^1$. The quotient surface  
$\Si:=X/\si$ is a normal regular surface with  $84k$ singular points of type $A_1$ that are
the images of the fixed points of $\iota\times Id$ on $X$. The
$\Z_3-$action on
$X$ commutes with
$\iota\times \Id$ and therefore there is an induced $\Z_3-$action on $\Si$. 
The fixed points of $\Z_3$ on $\Si$ are the $6$ image points of the fixed
points of $\Z_3$ on $X$ and are smooth for $\Si$. The automorphism  $\xi$ acts as the
multiplication by $\omega$  on the tangent  space at $3$ of these points  and it
has eigenvalues
$\omega$ and $\omega^2$ at the remaining
$3$.
If $T$ is the minimal desingularization of $\Si/\Z_3$,  then the usual computation gives:
$$K^2_T=24k-9;\quad q(T)=0;\quad p_g(T)=7k-1 .$$
Hence $p_g(Y)=p_g(T)$.

 One can use the same argument as in previous examples to show  that for $k\ge 1$ the
canonical
 map of $T$  is
 birational onto its image.
\bigskip
\begin{rem} Notice that no good canonical cover $X\to S$ such that $q(S)>0$ is known. It would
be interesting to see whether such an example exists.
\end{rem}

\section{The fundamental group}
In this section we prove that for the good canonical covers $Y\to T$ of Example 1 of \S 3 the
surface $Y$ is simply connected. We use the notation of \S 2 and \S 3.

Our main tool is the following:
\begin{thm}\label{fundgroup}
Let $Z$ be a connected simply connected  manifold, let $\Ga$ be a
group that acts properly discontinuously on $Z$. Assume that the stabilizer of every point
of $Z$ is finite and  let
$N<\Ga$ be the subgroup generated by the elements that do not act freely.
Then:

$$\pi_1(Z/\Ga)=\Ga/N.$$
\end{thm}
\begin{proof} The  quotient map $Z\to Z/\Ga$ has the path lifting property up to homotopy by
Proposition 3 of \cite{armstrong}. Hence the assumptions of Theorem 4 of \cite{armstrong} are
satisfied and we have the result.
\end{proof}

\begin{thm} Let $Y\to T$ be one of the good canonical covers of Example 1 of \S 3. Then:
$$\pi_1(Y)=0.$$
\end{thm}
\begin{proof}
We write $\bY:=X/\Z_3$, and we let $\eta\colon Y\to \bY$ be the minimal resolution of the
singularities of $\bY$.
The singularities of $\bY$ are points of type 
$A_1$,
$A_2$ or
$\frac{1}{3}(1,1)$. A standard application of van Kampen's theorem shows that
$\eta_*\colon\pi_1(Y)\to\pi_1(\bY)$ is an isomorphism. Hence from now on we study
$\pi_1(\bY)$.

 By the  Lefschetz Theorem
on hyperplane sections, the inclusion $X\to \Pic^1(C)\times \pp^1$ induces an isomorphism on
the fundamental groups. Therefore, if we denote by
$p\colon \C^2\times \pp^1\to \Pic^1(C)\times \pp^1$ the universal cover and we set
$\tX:=p\inv X$, then the restricted map $p\colon \tX\to X$ is the universal cover.
We identify the curve $C$ with $C\times (1:0)\subset\Pic^1(C)\times \pp^1$ and we set
$\tC:=p\inv C$. The restricted map $p\colon \tC\to C$ is the covering of $C$
associated to the commutator subgroup $[\pi_1(C),\pi_1(C)]$.
The automorphism $\xi$ of $\Pic^1(C)\times\pp^1$ lifts to the universal cover and we have an
extension of groups:
$$1\to\pi_1(\Pic^1(C)\times\pp^1)=\Z^4\to \Ga\to \Z_3\to 1$$
where $\Ga$ is a group that acts properly discontinuously  on $\C^2\times
\pp^1$. 
Clearly $\Ga$ acts also on $\tX$ and $\tC$ by construction, and $\tX/\Ga=\bY$,
$\tC/\Ga=C/\Z_3=\pp^1$.  If we denote by $N$ the subgroup of $\Ga$ generated by the elements
that do not act freely on $\tX$ then we have $\pi_1(\bY)=\Ga/N$ by Theorem \ref{fundgroup}.
We denote by $N_1$ the subgroup of $\Ga$ generated by the elements that do not
act freely on $\tC$. The fixed points of an element that does not act freely on
$\tC$ necessarily map to fixed points of $\Z_3$ on $C$, namely to one of the
points $Q_1\times (1:0)\dots Q_4\times (1:0)$. Since these four points are also
in $X$, it follows that $N_1\subset N$.
Consider now the universal cover $D\to C$ of $C$. The automorphism  $\xi$ of $C$
lifts to $D$ and we have an extension of groups analogous to the one above:
$$1\to\pi_1(C)\to \Ga_0\to\Z_3\to 1.$$
We denote by $N_0<\Ga_0$ the subgroup  generated by the elements that do not
act freely on $D$. Since $D/\Ga_0=C/\Z_3=\pp^1$ is simply connected, we have
$N_0=\Ga_0$ by  Theorem
\ref{fundgroup}. The commutator subgroup $[\pi_1(C),\pi_1(C)]$ is a
characteristic subgroup of $\pi_1(C)$, hence it is normal in $\Ga_0$. It
follows that every  element of
$\Ga_0$ induces an automorphism of $\tC=D/[\pi_1(C),\pi_1(C)]$.  Hence there is
a natural surjection 
$\Ga_0\to \Ga$.   Since the image of
$N_0$  in
$\Ga$ is obviously contained in $N_1$, we have $N_1=\Ga$. Since $N_1\subset
N$, we also have
$N=\Ga$, hence $Y$ is simply connected by Theorem \ref{fundgroup}.
\end{proof}

\bigskip

\noindent\begin{minipage}{13cm}
\parbox[t]{6.5cm}{Ciro Ciliberto\\
Dipartimento di Matematica\\
Universit\`a di Roma Tor Vergata\\
Via della Ricerca Scientifica\\
00133 Roma, Italy\\
cilibert@mat.uniroma2.it}
\hfill
\parbox[t]{5.5cm}{Rita Pardini\\
Dipartimento di Matematica\\
Universit\`a di Pisa\\
Via Buonarroti, 2\\
56127 Pisa, ITALY\\
pardini@dm.unipi.it}
\bigskip
\bigskip

\parbox[t]{6.5cm}{Francesca Tovena\\
Dipartimento di Matematica\\
Universit\`a di Roma Tor Vergata\\
Via della Ricerca Scientifica\\
00133 Roma, Italy\\
tovena@mat.uniroma2.it}
\hfill

\end{minipage}

\begin{thebibliography}{ABCD}
\bibitem[Ar]{armstrong} M. A. Armstrong, {\em Calculating the fundamental group of an orbit
space}, Proc. A.M.S. {\bf 84}, n. 2 (1982), 267--271. 

\bibitem[Ba] {bab} D. W. Babbage, {\em Multiple canonical surfaces}, Proc.
Cambridge
Phil. Soc., {\bf 30} (1934).

\bibitem[Be] {beauville} A. Beauville, {\em L'application canonique pour le
surfaces de type
g\'e\-n\'e\-ral}, Invent. Math., {\bf 55} (1979), 121--140.
\bibitem[Ca1]{babbage} F. Catanese, {\em Babbage's conjecture,  contact of
surfaces,
symmetric determinantal varieties and applications},  Invent. Math.,
{\bf 63} (1981),
433--465.
\bibitem[Ca2]{Ca} F. Catanese, {\em Canonical rings and ``special''
surfaces of general type},  Proceedings of Symposia in Pure Mathematics,
{\bf 46}
(1987), 175--194.
\bibitem[Ci] {enriques} C. Ciliberto, {\em A few comments on some aspects of
the
mathematical work of F. Enriques}, in Geometry and Complex Variables (S.
Cohen ed.),
Lecture Notes in Pure and Appl. Math, Marcel Decker (1991), 89-109.


\bibitem[CPT]{cirifra} C. Ciliberto, R. Pardini, F. Tovena, {\em Prym
varieties and the canonical map of surfaces of general type},  Annali della Scuola Normale
Superiore di Pisa. Classe di Scienze, XXIX (2000), 905--938.

\bibitem[MP] {MP} M. Mendes Lopes, R. Pardini, {\em Irregular canonical double
surfaces},
Nagoya J. of Math., {\bf 152} (1998), 203-230.


\bibitem[Pa1]{ritaabel}
R. Pardini,
{\em Abelian covers of algebraic
varieties},
 J. reine angew. Math. {\bf 417} (1991), 191--213.

\bibitem[Pa2] {supcan} R. Pardini, {\it Canonical images of surfaces}, J.
reine angew.
Math., {\bf 417} (1991), 215--219.
\bibitem[St]{steenbrink}J. H. M. Steenbrink, 
{\em  Mixed Hodge structure on the vanishing cohomology},
  in      ``Real and Complex Singularities, Oslo 1976'' (P. Holm ed.),  Sijthoff \&
Noordhoff, The Netherlands (1977),  525--564.

\bibitem[VdGZ]{zagier} G. van der Geer, D. Zagier, {\em The Hilbert modular
group for the
field ${\bf Q}(\sqrt{13})$}, Invent. Math., {\bf 42} (1977), 93--134.


\end{thebibliography}
\end{document}